\documentclass[11pt]{article}
\usepackage[centertags]{amsmath}
\usepackage{amsmath, amsfonts, amssymb, latexsym, dsfont}
\begin{document}

%%% Abbreviation %%%
\newcommand{\qed}{\hfill\Box}
\newcommand{\bi}{\bigskip}
\newcommand{\sm}{\smallskip}
\newcommand{\w}{\sqrt}
\newcommand{\wh}{\widehat}
\newcommand{\wt}{\widetilde}
\newcommand{\ee}{\end{equation}}
\newcommand{\eea}{\end{eqnarray}}
\newcommand{\bean}{\begin{eqnarray*}}
\newcommand{\eean}{\end{eqnarray*}}
\newcommand{\suml}{\sum\limits}
\newcommand{\intl}{\int\limits}
\newif\ifpctex
\newcommand{\noi}{\noindent}
\newcommand{\vt}{\vartheta}
\newcommand{\vp}{\varphi}
\newcommand{\vr}{\varrho}
\newcommand{\ve}{\varepsilon}

%%%%%%%  Theorem environment  %%%%%%%
  \newtheorem{definition}{Definition}
  \newtheorem{prob}{Problem}
  \newtheorem{conjecture}{Conjecture}
  \newtheorem{quest}{Question}
  \newtheorem{theorem}{Theorem}
  \newtheorem{proposition}{Proposition}
  \newtheorem{lemma}{Lemma}
  \newtheorem{corollary}{Corollary}
  \newtheorem{assumption}{Assumption}
  \newtheorem{example}{Example}
  \newtheorem{remark}{Remark}
  \newtheorem{result}{Result}
%%%%%%%  Domains of numbers   %%%%%%%
  \newcommand{\R}{\mathbb{R}}
  \newcommand{\Z}{\mathbb{Z}}
  \newcommand{\N}{\mathbb{N}}
  \newcommand{\C}{\mathbb{C}}
  \newcommand{\CL}{\mathcal{L}}
  \newcommand{\CH}{\mathcal{H}}
  \newcommand{\CM}{\mathcal{M}}
  \newcommand{\CU}{\mathcal{U}}

\newcommand{\Rand}[1]{\marginpar{#1}}
\renewcommand{\Rand}[1]{}
\marginparwidth2.5cm
\newcommand{\be}[1]{\Rand{\vspace{0,6cm}\tt #1}\begin{equation}\label{#1}}
\newcommand{\bea}[1]{\Rand{\vspace{0,6cm}\tt #1}\begin{eqnarray}\label{#1}}
\newcommand{\beL}[2]{\Rand{\vspace{0,6cm}\tt #1}\begin{lemma}[#2]\label{#1}}
\newcommand{\beA}[2]{\Rand{\vspace{0,6cm}\tt #1}\begin{assumption}[#2]\label{#1}}
\newcommand{\beD}[2]{\Rand{\vspace{0,6cm}\tt #1}\begin{definition}[#2]\label{#1}}
\newcommand{\beT}[2]{\Rand{\vspace{0,6cm}\tt #1}\begin{theorem}[#2]\label{#1}}
\newcommand{\beP}[2]{\Rand{\vspace{0,6cm}\tt #1}\begin{proposition}[#2]\label{#1}}
\newcommand{\beC}[2]{\Rand{\vspace{0,6cm}\tt #1}\begin{conjecture}[#2]\label{#1}}
\newcommand{\beCor}[2]{\Rand{\vspace{0,6cm}\tt #1}\begin{corollary}[#2]\label{#1}}

%%% Convergences %%%
  \newcommand{\atp}[2]{\genfrac{}{}{0pt}{}{#1}{#2}}
  \newcommand{\llto}{{_{\displaystyle{\atp{\ll}{t\to\infty}}}}}
  \newcommand{\Tto}{{_{\displaystyle{\atp{\Longrightarrow}{t\to\infty}}}}}
  \newcommand{\tto}{{_{\displaystyle{\atp{\longrightarrow}{t\to\infty}}}}}
  \newcommand{\tTo}{{_{\displaystyle{\atp{\longrightarrow}{T\to\infty}}}}}
  \newcommand{\Tno}{{_{\displaystyle{\atp{\Longrightarrow}{n\to\infty}}}}}
  \newcommand{\TNo}{{_{\displaystyle{\atp{\Longrightarrow}{N\to\infty}}}}}
  \newcommand{\tno}{{_{\displaystyle{\atp{\longrightarrow}{n\to\infty}}}}}
  \newcommand{\tMo}{{_{\displaystyle{\atp{\longrightarrow}{M\to\infty}}}}}
  \newcommand{\Mto}{{_{\displaystyle{\atp{\longrightarrow}{M\to\infty}}}}}
  \newcommand{\Nto}{{_{\displaystyle{\atp{\longrightarrow}{N\to\infty}}}}}
  \newcommand{\tNo}{{_{\displaystyle{\atp{\longrightarrow}{N\to\infty}}}}}
  \newcommand{\trhoo}{{_{\displaystyle{\atp{\longrightarrow}{\rho\to\infty}}}}}
  \newcommand{\Ttto}{{_{\displaystyle{\atp{\Longrightarrow}{T,t\to\infty}}}}}
  \newcommand{\Trhoo}{{_{\displaystyle{\atp{\Longrightarrow}{\rho\to\infty}}}}}

\title{Law of large numbers for superdiffusions: the non-ergodic case}
\thispagestyle{empty}
\author{J\'{a}nos Engl\"{a}nder\\
Department of Statistics and Applied Probability\\
University of California, Santa Barbara\\
CA 93106-3110, USA.\\
Email:{\texttt{englander@pstat.ucsb.edu}}}
\date{\today}
\maketitle

\begin{abstract}
In \cite{EngWin06} the Law of Large Numbers for the local mass of
certain superdiffusions was proved under a spectral theoretical
assumption, which is equivalent to the ergodicity (positive
recurrence) of the motion component of an $H$-transformed (or
weighted) superprocess. In fact the assumption is also equivalent to
the property that the scaling for the expectation of the local mass
is pure exponential.

In this paper we go beyond ergodicity, that is we consider cases
when the scaling is not purely exponential. \emph{Inter alia}, we
prove the analog of the Watanabe-Biggins Law of Large Numbers for
super-Brownian motion (SBM).

We will also prove another Law of Large Numbers for a bounded set
moving with subcritical speed, provided the variance term decays
sufficiently fast.

Further illustrative examples, such as SBM with drift and
super-Ornstein-Uhlenbeck process, will be provided too.
\end{abstract}
%\tableofcontents
\vspace{1cm} \noindent
{\it MSC} 2000 {\it subject classifications}.\ 60J60, 60J80\\
{\it Key words and phrases}.\ super-Brownian motion,
super-Ornstein-Uhlenbeck process, superdiffusion, superprocess, Law
of Large Numbers, $H$-transform, weighted superprocess, scaling
limit, local extinction, local survival.
\section{Introduction and statement of results}
\subsection{Basic notions}
Let $D\subseteq \mathbb{R}^d$ be a domain and let ${\cal B}(D )$
denote the  Borel sets of $D$. We write ${\cal M}_f(D)$ and ${\cal
M}_c(D)$ for the class of finite measures resp.\ the class of finite
measures with compact support on ${\cal B}(D )$. For $\mu\in{\cal
M}_f(D )$, denote $ \|\mu\|:=\mu (D )$ and let $C^+_b(D )$ and
 $C^+_c(D )$ be the class of non-negative bounded continuous
resp.\ non-negative continuous  functions $D \rightarrow\R$ having
compact support. Write $C^{k,\eta}(D )$ for the usual H\"older
spaces of index $\eta\in (0,1]$ including derivatives of order $k$,
and set $C^{\eta}(D):=C^{0,\eta}(D)$. Let $L$ be an elliptic
operator on the domain $D\subseteq\R^d$ of the form \be{A3}
   L:=\frac{1}{2}\nabla \cdot a \nabla +b\cdot\nabla,
\end{equation}
where $a_{i,j},b_i\in C^{1,\eta}(D)$, $i,j=1,...,d$, for some
$\eta\in (0,1]$, and the matrix $a(x):=(a_{i,j}(x))$ is symmetric,
and positive definite for all $x\in D$. In addition, let
$\alpha,\beta\in C^{\eta}(D)$, and assume that $\alpha$ is positive,
and $\beta$ is bounded from above.

\begin{definition}[$(L,\beta,\alpha;D )$-superdiffusion]\rm
With $D,L,\beta$ and $\alpha$ as above, let $\left(
X,\mathbf{P}^{\mu\,},\,\mu\in\mathcal{M}_f(D)\right) $ denote the
$(L,\beta,\alpha;D )$-superdiffusion. That is, $X$ is the unique
$\mathcal{M}_f(D)$-valued continuous (time-homogeneous) Markov
process which satisfies, for any $g\in C^+_b(D)$ ,
\begin{equation}
{\bf E}^{\mu}\exp\left\langle  X_{t\,},-g\right\rangle
 =\exp\, \langle\mu,-u(\cdot ,t)\rangle,
\label{Laplace.functional}
\end{equation}
where $u$ is the minimal nonnegative solution to
\begin{equation}
\left.
\begin{array}
[c]{c}%
u_{t}=Lu+\beta u-\alpha u^{2}\quad\text{\emph{on}}\;D%
\times(0,\infty),\vspace{8pt}\\
\lim\limits_{t\downarrow 0}u(\cdot,t)=g(\cdot).
\end{array}
\;\right\}  \label{cum.equ}%
\end{equation}
As usual, $\left\langle \nu,f\right\rangle $ denotes the integral
$\int_{D}\nu(\mathrm{d} x)\,f(x).$
\end{definition}

(See Dynkin (1991, 2002) \cite{Dyn91}, \cite{Dy02} or Dawson (1993)
\cite{Daw93} for the definition of superprocesses in general; see
Engl\"ander and Pinsky (1999) \cite{EngPin99} for more on the
definition in the particular setting above.)

One usually refers to $\beta$ as \emph{mass creation} and $\alpha$
as the \emph{intensity parameter} (or variance).

Note that although the above model is a time-homogeneous process,
later on the introduction of certain {\it time-inhomogeneous}
superdiffusions will be necessary.

\subsection{Law of Large Numbers vs. Local extinction}
Our principal interest is in establishing a Law of Large Numbers for
the local mass of certain superdiffusions. That is if $f\in
C_c^+(D)\not\equiv 0$ and $\mu\in \mathcal{M}_c(D)\not\equiv 0,$
then we would like to verify that,
$$\lim_{t\to\infty}\frac{\langle X_t,f \rangle}
{\mathbf{E}^{\mu}\langle X_t,f \rangle}= N_ {\mu},$$ where $N_
{\mu}$ is a (non-degenerate) random variable and the limit holds in
some suitable sense. (See \cite{EngWin06} for an explanation on why
such a statement can be called a Law of Large Numbers.)

One can immediately see however, that such a law cannot be true
without having some condition on the operator corresponding to the
superprocess. To elucidate on this point, note that  the Law of
Large Numbers will obviously not hold when $X$ {\em exhibits local
extinction} (i.e., the support of $X$ leaves any given bounded set,
${\bf P}^{\mu}$-a.s. for each $\mu\in{\cal M}_c$). Now, let \be{A5q}
\begin{aligned}
  \lambda _{c}=\lambda _{c}(L+\beta )
= \inf \{\lambda \in \mathbb{R}\ : \
  \exists u>0\ \text{satisfying}\ (L+\beta -\lambda )u=0\ \text{in}\
  D \}
\end{aligned}
\end{equation}
denote the {\it generalized principal eigenvalue }for $L+\beta $ on
$D$ (see the Appendix).  Since Pinsky proved that $X$ exhibits local
extinction if and only if $\lambda_c \le 0$, we can only hope to
have the Law of Large Numbers if we assume that $\lambda_c>0$.

\subsection{Motivation}
In order to understand what follows, we need to present first some
concepts regarding the \emph{criticality theory of second order
elliptic operators} (for a complete presentation, the reader is
referred to Chapter 4 in \cite{Pin95}). The operator
$L+\beta-\lambda_c$ is called \emph{critical} if there exists a
positive function $f$ satisfying that $(L+\beta-\lambda_c)f=0$ but
there is no (minimal positive) Green's function for the operator
$L+\beta-\lambda_c$. In this case $f$ is unique up to constant
multiples and is called the {\em ground state}. The operator
$L+\beta-\lambda_c$ is called product-critical if it is critical
with  ground state $0<\phi_c$, and $\phi_c$ and $\widetilde{\phi_c}$
(i.e.\ the ground state for the formal adjoint of
$L+\beta-\lambda_c$) satisfy $\langle
\mathrm{d}x,\phi_c\widetilde{\phi_c}\rangle<\infty$. In this case we
normalize them by $\langle
\mathrm{d}x,\phi_c\widetilde{\phi_c}\rangle =1$. If
$L+\beta-\lambda_c$ possesses a Green's function, then it is called
{\it subcritical}.
\smallskip

The following has been shown in \cite{EngWin06}:
\beP{JC1}{\cite{EngWin06}, Theorem 1} \label{LLNproduct-critical}

In addition to the assumption $\lambda_c>0$, also assume that
$L+\beta-\lambda_c$ is product-critical, that $\alpha\phi_c$ is
bounded and that $X$ starts in a state $\mu$ with
$\langle\mu,\phi_c\rangle<\infty$. Let $f\in C_c^+(D)$. If
$f\not\equiv 0$ and $||\mu||\neq 0$, then \be{A71}
\lim_{t\to\infty}\, \frac{\langle X_t,f\rangle}{{\bf E}^{\mu}\langle
X_t,f\rangle}
=\frac{\overline{X}^H_{\infty}}{\langle\mu,\phi_c\rangle}, \qquad
\mathrm{in }\;\,{\bf P^{\mu}}\mathrm{-probability}.
\end{equation}
\end{proposition}

A simple case of a superdiffusion is when $D=\mathbb R^d,\ d\ge 1,\
L=\frac{1}{2}\Delta$, with $\alpha,\beta$ positive constants
(supercritical super-Brownian motion). Here $\lambda_c=\beta$ and
$$\frac{1}{2}\Delta+\beta-\lambda_c= \frac{1}{2}\Delta.$$ Since
$\phi_c=\widetilde \phi_c\equiv 1,\,d\ge 1,$ the operator
$\frac{1}{2}\Delta$ is either critical but not product-critical
($d\le 2$), or subcritical ($d\ge 3$). Therefore \emph{this case was
not included in the setup} of \cite{EngWin06}. On the other hand,
the corresponding (Strong) Law of large Numbers is well known for
\emph{discrete particle systems}.

Using techniques from Fourier transform theory, Watanabe
\cite{Wat67} proved SLLN for branching-Brownian motion in $\mathbb
R^d$ and in certain subdomains of it. It is not clear however if his
method can be generalized for more general branching diffusions.
Furthermore, the proof in \cite{Wat67}  is thought to have a gap. In
\cite{Wat67} a family of nonnegative martingales
$\{W_t^{\lambda};\lambda\in\Lambda\}$ is considered, together with
the family of their limits
$\mathfrak{L}:=\{W^{\lambda};\lambda\in\Lambda\}$, where
$W^{\lambda}:=\lim_{t\to\infty}W_t^{\lambda}$. The problem, however,
is that those limits are only \emph{almost sure} ones, so for every
$\lambda\in\Lambda$ there is an exceptional null set $N_{\lambda}$.
Since $\Lambda$ is not countable, $\mathfrak{L}$ is not defined on
the \emph{uncountable} union $\bigcup_{\lambda\in\Lambda}
N_{\lambda}$!

Biggins \cite{Big92} gives a complete proof for the case when
branching-Brownian motion is replaced by branching RW. Here time is
discrete: $n=1,2,...$ and instead of considering $\{W_n^{\lambda}\}$
as a family indexed by $\lambda$, it is thought of as a sequence of
(continuous) functions $\{W_n(\cdot)\}$, i.e. $\lambda$ is not the
index but the argument. If the sequence is restricted to a compact
subset $F \subset \Lambda$, it can be thought of as a martingale
with values in the Banach space of continuous functions on $F$
(under the supremum norm). Biggins then proves that this martingale
converges almost surely and in mean.

The purpose of  this paper is to prove the Law of Large Numbers for
a class of superprocesses that \emph{includes supercritical
super-Brownian motion}. Instead of trying to adapt the
Watanabe-Biggins approach to our setting, our method will use some
ingredients from \cite{EngWin06}, however the time scales will have
to be modified now and also the spatial spread of the process must
be handled. In particular, to verify that supercritical
super-Brownian motion is included in the setup (satisfies the
assumption on the spatial speed), we will need a result from
\cite{Pin95b} too.

\subsection{Main results}
Throughout the paper the following assumption will be in force.
\begin{assumption}\rm
Let $\lambda_c=\lambda_c(L+\beta;D)$ denote the generalized
principal eigenvalue of $L+\beta$ on $D$, and let $S=\{S_t\}_{t\ge
0}$ denote the semigroup corresponding to $L+\beta-\lambda_c$ on
$D$.
\begin{enumerate}
\item [\textsf{(A.1)}] \textbf{(local survival)} $\lambda_c>0$,
\item [\textsf{(A.2)}] \textbf{(scaling of linear semigroup)} There exist two  functions $s:(0,\infty)\rightarrow(0,\infty)$,
and $h:D\rightarrow(0,\infty),\ h\in C^{2,\eta}(D)$,  and a locally
finite measure $r(\mathrm{d}x)$ such that
\begin{enumerate}
\item $\log s_t=\mathcal{O}(\log t)$ (i.e. $\sup_{t>0}\frac{\log s_t}{\log t}<\infty$)  and $\lim_{t\to\infty}(\log s_t)'=0,$
\item $\sup_D \alpha
h<\infty,$
\item $ \lim_{t\to\infty}\langle s_t\cdot S_t(f)(x),\mu (\mathrm{d}x)\rangle=
\langle f,r(\mathrm{d}y)\rangle\cdot \langle h,\mu \rangle$ for all
$f\in C_c^+(D)$, and $\mu\in \frac{1}{h}\mathcal{M}_f$,

\end{enumerate}

\item [\textsf{(A.3)}] \textbf{(spatial spread)} There exist two more functions,
$z,\wh z:(0,\infty)\rightarrow(0,\infty)$ such that
\begin{enumerate}
\item $\lim_{t\to\infty}\wh z_t=\infty$,
\item $\log (t+\wh z_t)=o(t),\ t\to\infty$,
\item $\lim_{t\to\infty} s_{t+\wh z_t}/s_{\wh z_t}=1$,
\item $\lim_{t\to\infty}P^{\mu}(\langle X_t,B^c_{z_t}(0)\rangle>0)=0$,
\item if
$f\in C_c^+(D)$, then $$\lim_{t\to\infty}\sup_{|x|\le z_t}\left|
\frac{s_{\wh {z_t}}}{h(x)}\cdot S_{\wh {z_t}}(f)(x) -\langle
f,r(\mathrm{d}y)\rangle\right|=0.$$

\end{enumerate}
\end{enumerate}
\end{assumption}
Let $H(x,t):=\exp (-\lambda_c t) h(x)$ and consider the weighted
superprocess $X^H$ (see section \ref{H-tr}). Abbreviate $W:={X}^H$
and let $\overline{W}$ denote the {\em total mass process}, i.e.
\be{Y11}
   \overline{W} =\overline{X}^H:=\| X^H\|.
\end{equation}
It is not hard to show that $\overline{W}$ is a supermartingale (see
the Appendix), and therefore it has a limit,
$\overline{W}_{\infty}$. One does not have $\mathbf{E}^{\mu}
\overline{W}_{\infty}=\langle \mu, h \rangle$ in general. However,
when $L^h_0$ is conservative on $D$, that is, it never leaves the
domain $D$ with probability one, $\overline{W}$ is a UI martingale
(see again the Appendix). In this case, by uniform integrability,
$\mathbf{E}^{\mu} \overline{W}_{\infty}=\langle \mu, h \rangle>0$.

\begin{theorem}[Law of Large Numbers]\label{LLN}
With the notations of Assumption 1, if $f\in C_c^+(D)\not\equiv 0$
and $\mu\in \mathcal{M}_c(D)\not\equiv 0,$ then in probability,
$$\lim_{t\to\infty}\frac{\langle X_t,f \rangle}
{\mathbf{E}^{\mu}\langle X_t,f \rangle}=
\frac{\overline{W}_{\infty}}{\langle \mu, h \rangle}.$$ The limit is
mean-one (and in particular, not identically zero) when $L_0^h$
corresponds to a conservative diffusion.
\end{theorem}
In order to give a simple condition for the limit to be mean-one, we
recall the \emph{compact support property}.

\begin{definition}
The $(L,\beta,\alpha;D )$-superdiffusion possesses the compact
support property if
\begin{equation}
{\mathbf P}^\mu\left(\bigcup_{0\le s\le t} \text{supp}\ \  X_s\
\subset\subset D\right)=1, \ \text{for all}\ \mu\in
\mathcal{M}_c(D),\ t\ge0.
\end{equation}

\end{definition}
[Here $A\subset \subset D$ means that the closure of the bounded
domain $A$  is  in $D$.] Since there are various conditions given in
\cite{EngPin99,EngPin06} for the compact support property to hold,
the following result is useful.
\begin{theorem}[No loss of mass in the limit]\label{ndl}
If the compact support property holds, then the  diffusion process
corresponding to $L_0^h$ on $D$ is  conservative, and consequently,
the limit appearing in Theorem \ref{LLN} is mean-one.
\end{theorem}

Using $H$-transforms  we will also prove another Law of Large
Numbers for a bounded set moving with subcritical speed, provided
the variance term $\alpha$ decays sufficiently fast.

\begin{theorem}[Law of Large Numbers for a moving bounded
set]\label{movingLLN} Let $c<\sqrt{2\beta}$, $\alpha(x)\le K
\exp(-cx)$ with some $K>0$, and define
$f^{(ct)}(x):=f(x_1+ct,x_2,...,x_d)$. With the notations of
Assumption 1, if $f\in C_c^+(D)\not\equiv 0$ and $\mu\in
\mathcal{M}_c(D)\not\equiv 0,$  then in probability,
$$\lim_{t\to\infty}\frac{\langle X_t,f^{(ct)} \rangle}
{\mathbf{E}^{\mu}\langle X_t,f^{(ct)}\rangle}= N_c,$$ where $N_c$ is
a nonnegative, mean one random variable.
\end{theorem}
\subsection{Outline}
The rest of this paper is organized as follows. In section \ref{Ex}
we give examples for our main theorem. In section \ref{H-tr} we give
a review on a fundamental tool, a space-time transformation. In
section \ref{Pr} we present the proofs. Finally, the Appendix
summarizes the necessary background.
\section{Examples}\label{Ex}
In this section we give five examples. All of them satisfy the
assumption, and thus obey the Law of Large Numbers. In all the
examples, $D=\mathbb R^d$, and the constant $c$ (appearing in all
but the first example) is positive. In our notation $x$  is the
$d$-dimensional vector $(x_1,x_2,...,x_d)$ and
$x^2:=|x|^2=\sum_{i=1}^{d} x_i^2$. All the examples are versions of
either the super-Brownian motion (SBM) or the
super-Ornstein-Uhlenbeck process (SOU).
\begin{example}[supercritical SBM]\rm The assumptions are satisfied for
supercritical super-Brownian motion. Indeed, if $\beta(\cdot)\equiv
\beta>0$, then  $\lambda_c=\beta$, because $\lambda_c(\Delta,\mathbb
R^d)=0$. Furthermore choose $h\equiv 1$,
$r(\mathrm{d}x):=\mathrm{d}x$ and the Brownian scaling factor
$s_t:=t^{d/2}$. Finally, as far as the spatial spread of the process
is concerned, $z_t$ can be defined as
$z_t:=(\sqrt{2\beta}+\epsilon)t$ (see \cite{Pin95b}), and thus $\wh
z_t$ can be defined e.g. as $\wh z_t=t^m$ with $m>2$. This setting
satisfies conditions (2)-(3), as long as $0<\alpha$ is bounded from
above.$\hfill \diamond$

Consider now the simplest case of the previous example, the one when
$\alpha$ is a positive constant. Then the non-degenerate random
variable $\overline{W}_{\infty}$ can be thought of as the scaled
limit of a one dimensional diffusion. Indeed, $Y:=\|X\|$ is a
diffusion corresponding to the operator
$$\mathcal{L}:=x\left(\alpha \frac{\partial^2}{\partial x^2}+\beta
\frac{\partial}{\partial x}\right)\ \mathrm{on}\ [0,\infty)$$ with
$Y_0=\|\mu\|$, and $\overline{W}_{\infty}=\lim_{t\to\infty}e^{-\beta
t}Y_t$.
\end{example}
For the reader unfamiliar with \emph{spatial $h$-transform} it is
helpful to review section \ref{H-tr} before reading the rest of the
examples. (The $h$-transform is a particular case of the
$H$-transform with $H(t,x)=h(x)$.)

In the following examples there is no need  to check our assumptions
for the Law of Large Numbers to hold. The validity of the Law of
Large Numbers will simply follow from its invariance under
$h$-transforms and from the Law of Large Numbers in the previous
example. This will in particular mean that the limiting random
variable is always \emph{non-degenerate}.

Let $X$ be the supercritical super-Brownian motion of the first
example. Since the Law of Large Numbers holds true for $X$ starting
with any nonzero finite measure, therefore it is also true for $X^h$
starting with any measure of the form $\nu=h\mu$, where
$0\not\equiv\mu$ is a finite measure. To avoid working with these
$h$-dependent spaces, we will simply assume in all the examples
below that the initial (nonzero) measure belongs to $\mathcal{M}_c$.
\begin{example}[supercritical SBM with drift]\rm
Let $h(x):=e^{cx_1}$. Then $X^h$ corresponds to the motion generator
$$\mathcal{L}:=\frac{1}{2}\Delta + b \cdot\nabla,\ \ b:=(c,0,...,0)$$
(Brownian motion with drift in the first coordinate direction), mass
creation $\beta^h=\beta+c^2/2>c^2/2$ and intensity parameter
$\alpha^h=\alpha e^{cx_1}$. Hence, \emph{the super-Brownian motion
with drift $c>0$ obeys the Law of Large Numbers, if the mass
creation is larger than $c^2/2$, and the intensity parameter is
$\mathcal{O}(e^{cx_1})$}, $|x_1|\to\infty$.$\hfill \diamond$
\end{example}
\begin{example}[supercritical SBM with outward drift]\rm
Let $0<h\in C^{2,\eta}$ satisfy $h(x):=e^{c|x|},\ |x|>>1$. Then
$X^h$ corresponds to the motion generator
$$\mathcal{L}:=\frac{1}{2}\Delta + b\cdot\nabla,\ \ b(x)=c\,
\frac{x}{|x|},\ \textrm{for} |x|>>1,$$ mass creation
$\beta^h=\beta+c^2/2>c^2/2$ and $\alpha^h=\alpha e^{c|x|}$. Hence,
\emph{any superprocess obeys the Law of Large Numbers, if the motion
component is Brownian motion with outward drift $c>0$ and the mass
creation is larger than $c^2/2$, while the intensity parameter is
$\mathcal{O}(e^{c|x|})$, $|x|\to\infty$.}$\hfill \diamond$
\end{example}
The last two examples concern `super-Ornstein-Uhlenbeck processes'.
\begin{example}[supercritical SOU with small $\alpha$]\rm
We now choose $h(x):=e^{-cx^2}$. Then $X^h$ corresponds to the
Ornstein-Uhlenbeck migration with generator
$$\mathcal{L}:=\frac{1}{2}\Delta - 2cx\cdot\nabla,$$  mass creation
$\beta(x)=K+2c^2x^2$, $K>-c$, and $\alpha^h:=\alpha e^{-cx^2}$.
Hence, \emph{the super-Ornstein-Uhlenbeck process with drift $c>0$
obeys the Law of Large Numbers, as long as the mass creation is of
the above form, and the intensity parameter is
$\mathcal{O}(e^{-cx^2})$, $|x|\to\infty$.}$\hfill \diamond$
\end{example}
\begin{example}[supercritical outward SOU with large $\beta$]\rm
Taking $h(x):=e^{cx^2}$, $X^h$ corresponds to the migration with
generator
$$\mathcal{L}:=\frac{1}{2}\Delta + 2cx\cdot\nabla,$$  mass creation
$\beta(x)=K+2c^2x^2$, $K>c$, and $\alpha^h=\alpha e^{cx^2}$. Hence,
\emph{any superprocess obeys the Law of Large Numbers, if the motion
component is an outward drifting Ornstein-Uhlenbeck process with
drift $c>0$ and the mass creation is of the above form, while the
intensity parameter is $\mathcal{O}(e^{cx^2})$,
$|x|\to\infty$.}$\hfill \diamond$
\end{example}
We note that if in the last two example one replaces $\beta(\cdot)$
by a positive constant $\beta$, then the models will belong to the
setup treated in \cite{EngWin06} (product-critical, or ergodic
case).
\section{The $H$-transform of superdiffusions}\label{H-tr}
\label{Htransformsection}\rm This section gives a review on the
$H$-transform\footnote{The reader should not confuse with the
space-time harmonic transformation yielding a Girsanov-type change
of measure -- see e.g. \cite{Ove94}.}. The $H$-transform, introduced
in \cite{EngWin06}, is a space-time generalization of the spatial
$h$-transform for superdiffusions ($h$-transformed, or weighted
superprocesses had been introduced earlier in \cite{EngPin99}).

We first review the more general definition of a  time-inhomogeneous
superdiffusion. Let $\widetilde{L}$ be an elliptic operator on
$D\times\R^+$ of the form \be{A3n}
   \widetilde{L}:=\frac{1}{2}\nabla \cdot \widetilde{a} \nabla +\widetilde{b}\cdot\nabla
\end{equation}
where the functions
$\widetilde{a}_{i,j},\widetilde{b}_i:D\times\R^+\to\R$,
$i,j=1,...,d$ are $C^{1,\eta}(D)$  (for some $\eta\in (0,1]$) in the
space, and continuously differentiable in the time coordinate.
Moreover assume that the symmetric matrix
$\widetilde{a}(x,t):=(a_{i,j}(x,t))$ is positive definite for all
$x\in D$ and $t\in\R^+$.

In addition, let
$\widetilde{\alpha},\widetilde{\beta}:D\times\R^+\to\R$, be
$C^{\eta}(D)$ in the space, and continuously differentiable in the
time coordinate. Finally assume that $\widetilde{\alpha}$ is
positive, and $\widetilde{\beta}$ is bounded from above.

\beD{D2}{Time-inhomogeneous
$(\widetilde{L},\widetilde{\beta},\widetilde{\alpha};D)$-superdiffusion}
\rm\
\begin{itemize}

\item[{}] {\bf (i)} The
{\rm$(\widetilde{L},\widetilde{\beta},\widetilde{\alpha};D)$-superdiffusion}
is a  measure-valued (inhomogeneous) Markov process, $(X,{\bf
P}^{\mu,r};\,\mu\in{\cal M}_f(D),r\ge 0)$,  that is, a family
$\{{\bf P}^{\mu,r}\}$ of probability measures where ${\bf
P}^{\mu,r}$ is a probability on $C([r,\infty) )$ and the family is
indexed by ${\cal M}_f(D))\times [0,\infty)$, such that the
following holds: for each $g\in C_b^+(D)$ and $\mu\in{\cal M}_f(D)$,
\be{A4c} {\bf E}^{\mu,r}[\exp-\langle X_t,g\rangle]=
   \exp-\langle\mu,u(\cdot,r;t,g)\rangle,
\end{equation}
where $u=u(\cdot,\cdot;t,g)$ is a particular  non-negative solution
to the backward equation \be{A5s}
\begin{aligned}
   -\partial_r u&=\widetilde{L}u+\widetilde{\beta} u-\widetilde{\alpha} u^2\hspace{1cm}
  \mbox{in }D\times (0,t),\\
   \lim_{r\uparrow t}u(\cdot,r;t,g)&=g(\cdot).
\end{aligned}
\end{equation}
{\bf (ii)} To determine the solution $u$ uniquely, use the
equivalent {\it forward} equation along with the minimality of the
solution: fix $t>0$ and introduce the `time-reversed' operator
$\wh{L}$ on $D\times (0,t)$ by \be{A3n2}
   \wh{L}:=\frac{1}{2}\nabla \cdot \hat{a} \nabla
   +\hat{b}\cdot\nabla,
\end{equation}
where, for $r\in[0,t]$,
$$\hat a(\cdot,r):=\widetilde a(\cdot,t-r)\ \mathrm{and}\ \hat b(\cdot,r):=\widetilde b(\cdot,t-r);$$
furthermore let
$$\hat \beta(\cdot,r):=\widetilde \beta(\cdot,t-r)\ \mathrm{and}\ \hat \alpha(\cdot,r):=
\widetilde \alpha(\cdot,t-r).$$ Consider now $v$, the {\it minimal}
non-negative solution to the {\it forward} equation\be{A5s2}
\begin{aligned}
   \partial_r v&=\wh{L}v+\hat{\beta} v-\hat{\alpha} v^2\hspace{1cm}
  \mbox{in }D\times (0,t),\\
   \lim_{r\downarrow 0}v(\cdot,r;t,g)&=g(\cdot).
\end{aligned}
\end{equation}
Then $$u(\cdot,r;t,g)=v(\cdot,t-r;t,g).$$ (See also \cite{EngWin06}
concerning the construction of minimal non-negative solutions for
forward equations.)
\end{itemize}
\end{definition}

As we will see in Lemma \ref{H-transform} (b), one way of defining a
time-inhomogeneous superdiffusion is to start with a
time-homogeneous one, and then to apply an `{\it $H$-transform}'. In
general, the $H$-transform of a time-inhomogeneous superdiffusion is
defined as follows. Let $0<H\in C^{2,\eta}(D)\times
C^{1,\eta}(\mathbb R^+)$ and let $X$ be a
$(\widetilde{L},\widetilde{\beta},\widetilde{\alpha};D)$-superdiffusion.
We define a new process $X^H$ by
\begin{equation}\label{newprocess}
X_t^H:=H(\cdot,t)X_t\quad \left(\mbox{that is,}\
\frac{\mathrm{d}X_t^H}{\mathrm{d}X_t}=H(\cdot,t)\right), \quad t\geq
0.
\end{equation}
This way one obtains a new superdiffusion, which, in general, {\it
is not finite measure-valued} but only $\sigma$-finite
measure-valued. That is, if $\mathcal{M}(D)$ denotes the family of
all (finite or infinite) measures on $D$, then
$$X^H_t\in\mathcal{M}_H^{(t)}(D):=\{\nu\in \mathcal{M}(D)\mid H(\cdot,t)^{-1}\nu\in
\mathcal{M}_f(D)\}$$ (c.f. \cite{EngPin99}, p. 688.)

In \cite{EngPin99}, Section 2, it was shown, that, from an
analytical point of view, the (spatial) $h$-transform of the
superdiffusion is given by a certain transformation of the
corresponding semilinear operator. This remains the case for the
space-time $H$-transform. The following result was proved in
\cite{EngWin06}.
\begin{lemma}\label{H-transform}
\begin{itemize}
\item[{}] \hspace{2cm}

Let $X^H$ be defined by (\ref{newprocess}). Then \item[(a)] $X^H$ is
a $\left(\widetilde{L}+\widetilde{a}\frac{\nabla H}{H}\cdot\nabla,
  \widetilde{\beta}+\frac{\widetilde{L}H}{H}+\frac{\partial_rH}{H},
  \widetilde{\alpha}H;D\right)$-superdiffusion.
\item[(b)] In particular, if $X$ is a time-homogeneous
  $(L,\beta,\alpha;D)$-superdiffusion, and $H$ is of the form
\be{H}
   H(x,t):=e^{-\lambda_c t}h(x),
\end{equation}
where $\lambda_c$ is the principal eigenvalue of $L+\beta$, and $h$
is a positive solution of $(L+\beta)h=\lambda_c h$, then $X^H$ is a
$(L+a\frac{\nabla h}{h}\cdot \nabla,0,\alpha h e^{-\lambda_c
t};D)$-superdiffusion.
\end{itemize}
\end{lemma}
\begin{remark}[Unbounded $\widetilde \beta$'s]
\rm As it is already the case with the spatial $h$-transform for
superdiffusions, it is possible that the coefficient $\widetilde
\beta$ transforms into a new coefficient that is no longer bounded.
In fact this can be the very definition of superdiffusions for
certain unbounded $\widetilde \beta$'s (see \cite{EngPin99}, Section
2 for explanation).$\hfill \diamond$
\end{remark}
\begin{remark}[Invariance of LLN]\rm
When the function $H$ is of the form $H(x,t)=h(x)r(t)$, the validity
of the Law of Large Numbers is \emph{invariant} under the
transformation. Indeed, if $f$ belongs to $C^+_c(D)$, then so does
$g:=hf$; and for $\nu=r(0)h\mu$,
$$\frac{\langle X^H_t,f \rangle}
{\mathbf{\widetilde{E}}^{\nu}\langle X^H_t,f \rangle}=\frac{\langle
X_t,g \rangle} {\mathbf{E}^{\mu}\langle X_t,g \rangle},$$ and
 the invariance follows by letting $t\to\infty$. $\hfill
\diamond$\end{remark}
\section{Proofs}\label{Pr}
\subsection{Proof of Theorem\ref{LLN}} Having the $H$-transform at hand we
now reformulate our assumptions and also the statement in terms of
the space-time weighted superprocess $X^H$, where $H$ is of the form
\be{Hforus}
   H(x,t):=e^{-\lambda_c t}h(x),
\end{equation} $\lambda_c$ is the principal eigenvalue of $L+\beta$,
and $h$ is the function appearing in the assumption. Abbreviate
\be{LH} L^{h}_0:=L+a\frac{\nabla h}{h}\cdot \nabla
\end{equation}
and note that in fact
$$L^{h}_0(u)=h^{-1}(L+\beta-\lambda_c)(h u)=
   H^{-1}(L+\beta+\partial_t)(Hu).$$

Let $\mathfrak{S}=\{\mathfrak{S}_s\}_{s\ge 0}$ denote the semigroup
corresponding to the operator $L_0^{h}=(L+\beta-\lambda_c)^h$, that
is, $\mathfrak{S}:=S^{h}$. Since the semigroup $\mathfrak{S}$
corresponds to an operator  with no zeroth order part (i.e. to
$L_0^{h}$),
\begin{equation}\label{contraction} \mathfrak{S}_s
1\le 1.
\end{equation}
If the diffusion process corresponding to $L^{h}_0$ on $D$ is
conservative,  then in fact $\mathfrak{S}_s 1= 1;$ in general one
only has (\ref{contraction}). (See more on conservativeness and its
connection to a uniformly integrable martingale in the Appendix).

Let us use the shorthand $u':=\partial_t u$. We claim that $h$ is a
positive solution of $(L+\beta)h=\lambda_c h$. To see this note that
(A.2)(c) with $\mu=\delta_x$ yields $\lim_{t\to\infty} s_{ t}\cdot
S_t(f)(x)= \langle f,r(\mathrm{d}y)\rangle\cdot h(x)$ for $x\in D$
and $f\in C_c^+(D)$. Then, defining
$\widehat\beta(t,x):=\beta(x)+(\log s_t)'$, the equation
$(L+\beta)h=\lambda_c h$ follows from the fact that
$u(t,x):=S_t(f)(x)$ solves
$$(L+\beta-\lambda_c)u=u'$$ and therefore
$v(t,x):=s_t\cdot u(t,x)$ solves
$$(L+\beta-\lambda_c)v=v'-(\log s_t)'v,$$
that is
$$(L+\widehat\beta-\lambda_c)v=v'.$$
Indeed, first note that
$$\lim_{t\to\infty}\sup_{x\in D}\left|\widehat\beta(t,x)-\beta(x)\right|=0.$$ Then a standard
argument (see \cite{EngPin99}, p.708) together with the second
relation in (A.2)(a) gives that
$w(x):=\lim_{t\to\infty}v(t,x)=\langle f,r(\mathrm{d}y)\rangle\cdot
h(x)$ belongs to $C^{2,\eta}$ and solves the steady state equation
$(L+\beta-\lambda_c)w=0$. Then, also
\begin{equation}\label{steady.state}
(L+\beta-\lambda_c)h=0.
\end{equation}
Since (\ref{steady.state}) holds,  $X^H$ is a $(L+a\frac{\nabla
h}{h}\cdot \nabla,0,\alpha h e^{-\lambda_c t};D)$-superdiffusion by
Lemma \ref{H-transform}(b). From now on $\widetilde{\mathbf P}$ will
denote the probability corresponding to $X^H$. By recalling the
definition of $\mathfrak{S}$  and defining
$$\nu:=h\mu,\ g:=f/h,\ q:=hr,$$
one can reformulate (A.2)(c) as follows:
$$\textsf{({A*}.2)(c)} \lim_{t\to\infty}\langle s_t\cdot \mathfrak{S}_t(g)(x),\nu
(\mathrm{d}x)\rangle=\langle g,q(\mathrm{d}y)\rangle\cdot \|\nu\| ,
\hspace{0.2cm} g\in C_c^+(D),\ \nu\in \mathcal{M}_f.$$ Similarly,
(A.3)(d)-(A.3)(e) become
\begin{itemize}
\item[]\textsf{({A*}.3)(d)} $\lim_{t\to\infty}\widetilde P^{\nu}
(\langle W_t,B^c_{z_t}(0)\rangle>0)=0$, and
\item[]\textsf{({A*}.3)(e)}  $\lim_{t\to\infty}\sup_{|x|\le
z_t}\big|s_{\wh {z_t}}\cdot \mathfrak{S}_{\wh {z_t}}(g)(x) -\langle
g,q(\mathrm{d}y)\rangle\big|=0,$  $g\in C_c^+(D)$.

\end{itemize}
Finally, the theorem itself transforms into the following statement:
if $g\in C_c^+(D)\not\equiv 0$ and $\nu\in
\mathcal{M}_c(D)\not\equiv 0,$ then in probability,
$$\lim_{t\to\infty}\frac{\langle X^H_t,g \rangle}
{\mathbf {\widetilde E}^{\nu}\langle X^H_t,g \rangle }=
\frac{\overline{W}_{\infty}}{\|\nu\|},$$ or, equivalently, in
probability,
\begin{equation}\label{equivalent}\lim_{t\to\infty}s_t\langle
W_t,g \rangle= \langle g , q(\mathrm{d}y) \rangle
\overline{W}_{\infty}.
\end{equation}
(In the equivalence we used that $\mathfrak{S}$ is the
\emph{expectation} semigroup along with (A*.2)(c).)

In order to show (\ref{equivalent}), the main idea is to use the
comparison with the deterministic flow as in \cite{EngWin06},
nevertheless, there is an essential difference. In \cite{EngWin06}
we argued that by considering some large time $t+T$ (where both $t$
and $T$ are large), the changes of $\overline{X}^H$ are negligible
after time $t$, while  the remaining time  $T$ is still long enough
to distribute the produced mass according to the ergodic flow given
by the $H$-transformed migration. We then let first $T$ and then $t$
go to infinity.

Reading carefully the proof in \cite{EngWin06} one can see that this
method  relied heavily on the ergodicity of the flow and would break
down here. Hence, instead of letting first $T$ and then $t$ go to
infinity, we now define
$$T:=T_t=\wh z_t.$$
Similarly to \cite{EngWin06}, the strategy is to
\begin{itemize}
\item[(\textsf{a})] show that the total mass more or less stabilizes
by time $T_t$,
\item[(\textsf{b})] identify the limit of the \emph{scaled flow}
(starting from the state of the
process at $T_t$) at time $t+T_t$
\item[(\textsf{c})] show that it agrees with the scaling limit of
the measure-valued process itself.
\end{itemize}
Of course, (\textsf{a}) is simple:  being a supermartingale, the
\emph{total mass converges}:
\begin{equation}\label{t.mass.conv}
\lim_{t\to\infty}\|W_t\|=\overline{W}_{\infty},\ \
\widetilde{\mathbf P}^{\nu}-a.s.
\end{equation}
Unlike in \cite{EngWin06}, however, we do not know \emph{a priori},
that the limit is non-zero, and moreover, one cannot proceed further
without exploiting what is known about the radial speed of the
process. Therefore we continue as follows. Let $\{Z_{W_t}(s)\}_{s\ge
0}$ denote the deterministic flow starting from the (random) measure
$W_t$. Since given $W_t$, $$\langle Z_{W_t}(\wh z_t),g
\rangle=\langle\mathfrak{S}_{\wh z_t}(g)(x),
W_t(\mathrm{d}x)\rangle,$$  one has
\begin{equation}\label{linear.spread}
\widetilde{\mathbf P}^{\nu}\left(\left|s_{\wh z_t}\langle
Z_{W_t}(\wh z_t),g \rangle-\|W_t\|\langle g,
q(\mathrm{d}y)\rangle\right|>\epsilon\right)\le A_t+B_t,
\end{equation}
where
$$A_t:=\widetilde{\mathbf P}^{\nu}\left[\|W_t\|\cdot
\sup_{|x|\le z_t}\left|s_{\wh z_t}(\mathfrak{S}_{\wh
z_t}(g))(x)-\langle g, q(\mathrm{d}y)\rangle\right|
>\epsilon\right],$$
$$B_t:=\widetilde{\mathbf P}^{\nu}\left[W_t(B^c_{z_t}(0)) >0\right].$$
By (A*.3)(d-e), one has
$\lim_{t\to\infty}A_t=\lim_{t\to\infty}B_t=0$. Hence,
\begin{equation}\label{sum.to.zero}
\lim_{t\to\infty}\widetilde{\mathbf P}^{\nu}\left(\left|s_{\wh
z_t}\langle Z_{W_t}(\wh z_t),g \rangle-\|W_t\|\langle g,
q(\mathrm{d}y)\rangle\right|>\epsilon\right)=0.
\end{equation}
From this, along with  (A.3)(c) and (\ref{t.mass.conv}), one obtains
the \emph{scaling limit of the flow}:
\begin{equation}\label{sc.limit.flow}
\lim_{t\to\infty}\widetilde{\mathbf P}^{\nu}\left(\left|s_{t+\wh
z_t}\langle Z_{W_t}(\wh z_t),g \rangle-\overline{W}_{\infty}\langle
g, q(\mathrm{d}y)\rangle\right|>\epsilon\right)=0.
\end{equation}
Our goal is therefore to show that \emph{the scaling limit of the
flow agrees with the scaling limit of the measure-valued process}.
To achieve  this, recall (A.2)(b). A computation using Chebysev's
inequality and the supermartingale property (essentially the same
computation as the one giving formula (28) in \cite{EngWin06})
yields:
\begin{eqnarray*}\label{Chebysev}
\widetilde{\mathbf P}^{\nu}\left(\left|s_{t+\wh z_t}\langle
Z_{W_t}(\wh z_t),g \rangle-s_{t+\wh z_t}\langle W_{t+\wh
z_t,g}\rangle\right|>\epsilon\right)\le C\, \frac{\widetilde{\mathbf
E}^{\nu}\|W_t\|}{\epsilon^2 \lambda_ce^{\lambda_c t}}\cdot
s^2_{t+\wh z_t}\\ \le C\, \frac{\|\nu\|}{\epsilon^2
\lambda_ce^{\lambda_c t}}\cdot s^2_{t+\wh z_t},
\end{eqnarray*}
where $C=C(\|g\|,\|\alpha h\|):=18\,\| \alpha h\|\cdot\|g\|^2$.
Recall the abbreviation $T:=\wh z_t$. If we show now that
$$\lim_{t\to\infty}e^{-\lambda_c t}s^2_{t+T}=0,$$
then we are done, since this implies
\begin{eqnarray*}\label{the.end}
\lim_{t\to\infty}\widetilde{\mathbf
P}^{\nu}\left(\left|s_{t+T}\langle Z_{W_t}(T),g
\rangle-s_{t+T}\langle W_{t+T,g}\rangle\right|>\epsilon\right)=0.
\end{eqnarray*}

To estimate $e^{-\lambda_c t}s^2_{t+T}$, observe that its logarithm
can be estimated by the first condition in (A.2)(a):
$$\log \left(e^{-\lambda_c t}s^2_{t+T}\right)=
-\lambda_c t +2 \log s_{t+T}\le -\lambda_c t +2k \log (t+T),$$ with
some $k>0$ and, by (A.3)(b),
$$\lim_{t\to\infty}\left[ -\lambda_c t +2k \log (t+T)\right]=
\lim_{t\to\infty}t\left[ -\lambda_c  +2k \log (t+T)/t\right]=
-\infty.$$ This completes the proof. $\qed$
\subsection{Proof of Theorem\ref{ndl}}

In \cite{EngPin06} the following was shown.
\begin{proposition}[\cite{EngPin06}, Theorem 3]\label{explosionnocompact}
Assume that the  diffusion process corresponding to $L$ on $D$ is
not conservative and that
\begin{equation}
\sup_{x\in D}\alpha(x)<\infty\ \text{and}\ \inf_{x\in
D}\beta(x)>-\infty.\label{alphabeta}
\end{equation}
Then  the compact support property does not hold.
\end{proposition}
(In \cite{EngPin06} this was stated for $D=\mathbb R^d$, but the
proof goes through for general $D$ too.)

Using this result, Theorem \ref{ndl} follows easily by applying an
$h$-transform. Indeed, let us suppose that $X$ possesses the compact
support property but the diffusion corresponding to $L_0^h$ is not
conservative. Since the support of the superprocess (and thus the
compact support property too) is invariant under $h$-transforms,
therefore $X^h$ possesses the compact support property too. On the
other hand, since $X^h$ is the $(L_0^h,\lambda_c,\alpha
h;D)$-superprocess, $L_0^h$ is not conservative, and  $\alpha h$ is
bounded from above, Proposition \ref{explosionnocompact} implies
that the compact support property does not hold;
contradiction.$\qed$

\subsection{Proof of Theorem\ref{movingLLN}}
Recall that in our notation $x$  is the $d$-dimensional vector
$(x_1,x_2,...,x_d)$ and let $h(x):=\exp(cx_1)$. A straightforward
computation reveals that $X^h$ is the superprocess corresponding to
the operator
$$\left[\frac{1}{2}\Delta +c\frac{\partial}{\partial x_1} +
\left(\beta +\frac{c^2}{2}\right) \right]u-\alpha h u^2,$$ and
starting at $\nu:=h\mu$. Therefore, $X^h$ in a coordinate system
$$x_1'=x_1-ct,\ x_2'=x_2,...,\ x_d'=x_d$$ is equal in distribution to the
super-Brownian motion corresponding to the operator
$$\left[\frac{1}{2}\Delta  +
\left(\beta +\frac{c^2}{2}\right) \right]u-\alpha h u^2.$$ (This can
be derived easily for example from the log-Laplace equation
(\ref{Laplace.functional}).) By Theorem \ref{LLN}, this latter one
obeys the Law of Large Numbers because $\alpha(x) \exp(cx_1)$ is
bounded and $c<\sqrt{2\beta}$ guarantees that $\lambda_c>0$.  (Of
course the limiting random variable depends on $c$.) For $X^h$ this
means that if $f\in C_c^+(D)\not\equiv 0$ and $\mu\in
\mathcal{M}_c(D)\not\equiv 0,$ then in probability,
$$\lim_{t\to\infty}\frac{\langle X^h_t,f^{(ct)} \rangle}
{\mathbf{\widetilde E}^{\nu}\langle X^h_t,f^{(ct)}
\rangle}=\frac{\overline{W}_{\infty}}{\| \nu \|},$$ where
$$f^{(ct)}(x):=f(x_1+ct,x_2,...,x_d),$$ and
$\overline{W}_{\infty}=\overline{W}_{\infty,c}$. Denoting $g:=hf$,
we obtain that if $g\in C_c^+(D)\not\equiv 0$ and $\mu\in
\mathcal{M}_c(D)\not\equiv 0,$ then in probability,
$$\lim_{t\to\infty}\frac{\langle X_t,g^{(ct)} \rangle}
{\mathbf{E}^{\mu}\langle
X_t,g^{(ct)}\rangle}=\frac{\overline{W}_{\infty}}{\langle h, \mu
\rangle},
$$
where $g^{(ct)}(x):=g(x_1+ct,x_2,...,x_d)$.$\qed$

\section{Appendix}
In this section we provide some background material regarding
various probabilistic and analytic concepts. These can be found more
completely in \cite{EngPin99,EngWin06,Pin95}.
\subsection{The particle picture for the superprocess}
In the introduction we defined the $(L,\beta,\alpha;D)$-superprocess
$X$ analytically. In fact, $X$ also arises as the short life time
and high density diffusion limit of a \emph{branching particle
system}, which can be described as follows: in the $n^{\mathrm{th}}$
approximation step each particle has mass $1/n$ and lives a random
time which is exponential with mean $1/n$. While a particle is
alive, its motion is described by a diffusion process  corresponding
to the operator $L$. At the end of its life, the particle dies and
is replaced by a random number of particles situated at the parent
particle's final position. The distribution law of the number of
descendants is spatially varying such that the mean number of
descendants is $1+\frac{\beta(x)}{n}$, while the variance is assumed
to be $2\alpha(x)$. All these mechanisms are independent of each
other.

See Appendix A in \cite{EngPin99} for a precise statement on the
particle approximation.
\subsection{The generalized principal eigenvalue}
Let $Y$ be the diffusion process on $D$ corresponding to $L$, and
denote by $\mathbb P^x$ the law of $Y$ starting at $x\in D$. Then
from a probabilistic point of view, the generalized principal
eigenvalue can be equivalently expressed as \be{A5qq}
  \lambda_{c}
 =
  \sup_{\{A:\ A\subset \subset D,\ \partial A\ \mathrm{is}\ C^{2,\eta}\}}
  \lim_{t\to \infty }\frac{1}{t}\log{\mathbb{E}^{x}
  \left[
  \exp\left(\int_{0}^{t}\beta (Y_s)\,\mathrm{d}s\right);\,\tau^{A}>t
  \right]},
\end{equation}
for any $x\in D,$ where $\tau ^{A}:=\inf \{ t\geq 0:Y_t\not\in A\}$,
and the $C^{2,\eta}$-boundary is defined with the help of
$C^{2,\eta}$-maps in the usual way.  (See Section 4.4 in
\cite{Pin95} on the subject). Hence, since $\beta$ is bounded from
above, $\lambda_{c}<\infty$; and it is known from standard theory
that for any $\lambda \geq \lambda _{c}$, there exists a function
$0<f \in C^{2,\eta }(D)$ such that $(L+\beta )f =\lambda f $ on $D$.
(See Section 4.3 in \cite{Pin95}.)

\subsection{Uniform integrability and conservativeness}
Just like in \cite{EngWin06}, we show that
\begin{equation}\label{bdd.variance}
\lim_{t\to\infty}\widetilde{\mathrm{Var}}^{h\mu}
(\overline{X}^H)=\int_0^{\infty}\mathrm{d}s\,e^{-2\lambda_c s}\,
\langle \mu,{\cal S}_s[\alpha h^2]\rangle<\infty.
\end{equation}
We will also see that if the diffusion process $\widehat{Y}$
corresponding to $L^{h}_0$ on $D$ is conservative, then
$\overline{X}^H$ is a uniformly integrable $\widetilde{\bf
P}^{h\mu}$-martingale, whereas in general,  $\overline{X}^H$ is a
(non-negative) $\widetilde{{\bf P}}^{h\mu}$-supermartingale.
Clearly, only in the latter case it becomes an issue whether the
limit is non-zero.

Indeed, if $\widehat{Y}$ is conservative, then consider the class
$$\mathfrak{C}(D):=\{f\in C^2(D):\ \exists \mho\subset
D\mathrm{\ bounded\ s.t.}\ \overline{\mho}\subset D;\
f=\mathrm{const}\ \mathrm{on}\ D\setminus \mho \}.$$ By
Lemma~\ref{H-transform}(b) along with Theorem A2 in \cite{EngPin99},
we have that for all $f\in \mathfrak{C}(D)$, \be{H3a}
\mathrm{d}\langle X^H_t,f\rangle =\langle
X^H_t,L^{h}_0f\rangle\,\mathrm{d}t+\mathrm{d}M_t(f),
\end{equation}
where $\{M_t(f)\}_{t\ge 0}$ is a square-integrable $\widetilde{{\bf
P}}^{h\mu}$-martingale, and its quadratic variation  $\langle
M(f)\rangle$ is given by \be{H4}
   \langle M(f)\rangle_t=\int^t_0\mathrm{d}s\,e^{-\lambda_c s}
   \langle X^H_s,\alpha
{h}f^2\rangle,\ \ t\ge 0.
\end{equation}
(The point is that one can take the function class $\mathfrak{C}(D)$
instead of just $C^2_c(D)$ when $\widehat{Y}$ is conservative.)

Applying (\ref{H3a}) to the function $f\equiv 1$, it follows that
$\overline{X}^H$ is a $\widetilde{{\bf P}}^{h\mu}$-martingale.
Furthermore,  by (\ref{H4}), \be{H4f<}
%A%
  \widetilde{{\bf E}}^{h\mu}\left[\langle X^H_t,1\rangle^2\right]
 =
 \langle\mu, h\rangle^2+\int_0^{t}\mathrm{d}s\,e^{-\lambda_c s}\,
 \langle h\mu,\mathfrak{S}_s[\alpha h]\rangle.
\end{equation}
That is \be{H4f<2} \widetilde{\mathrm{Var}}^{h\mu}
(\overline{X}^H_t)
 =\int_0^{t}\mathrm{d}s\,e^{-\lambda_c s}\,
   \langle h\mu,\mathfrak{S}_s[\alpha h]\rangle=
   \int_0^{t}\mathrm{d}s\,e^{-2\lambda_c s}\,
 \langle \mu, S_s[\alpha h^2]\rangle.
\end{equation}
Letting $t\to\infty$ we obtain (\ref{bdd.variance}). Replacing $t$
by $\infty$ in the first of the integrals in (\ref{H4f<2}), we have
from (\ref{contraction}) and from our assumptions that
$$\widetilde{\mathrm{Var}}^{h\mu} (\overline{X}^H_t)\le
\int_0^{\infty}\mathrm{d}s\,e^{-\lambda_c s}\,
   \langle h\mu,\mathfrak{S}_s[\alpha h]\rangle
\le{\lambda_c }^{-1}\parallel\alpha h\parallel_{\infty}\,
   \langle \mu, h\rangle<\infty.$$
Hence, by (\ref{H4f<}), $\sup\nolimits_{t\ge 0}\widetilde{{\bf
E}}^{\mu h}(\overline{X}^H_t)^2 <\infty,$ and consequently
$\overline{X}^H$ is uniformly integrable.

Let $\mathfrak{D}:=D\cup \{\Delta\}$ be the \emph{one-point
compactification} of $D$ (when the underlying diffusion process
$\widehat{Y}$ is non-conservative on $D$, $\Delta$ is the
\emph{cemetery state} for $\widehat{Y}$). Relaxing the assumption on
the conservativeness of $\widehat{Y}$,  the argument in
\cite{EngPin99}, pp. 726--727 shows that, although one can
\emph{not} work directly with the function class $\mathfrak{C}(D)$
(only with its subclass $C^2_c(D)$),
 extending ${X}^H$ with $\widetilde{{\bf P}}$
 appropriately onto $\mathfrak{D}$ with ${\mathfrak{P}}$
 makes $\overline{X}^H$ a
${\mathfrak{P}}^{h\mu}$-martingale. Now, since the mass on the
cemetery state $\Delta$ is nondecreasing in time, therefore
$\overline{X}^H$ is less than this martingale by a non-decreasing
process, that is, it is a $\widetilde{{\bf
P}}^{h\mu}$-\emph{supermartingale}. (In the non-conservative case,
intuitively, mass is `lost' at the Euclidean boundary of $D$ or at
infinity.)

\bigskip {\bf Acknowledgement.} The author owes thanks to J. Biggins
for a helpful discussion regarding the proofs in \cite{Big92} and
\cite{Wat67}.
 %\small

\end{document}